\documentclass[a4paper,sort&compress]{amsart}
\usepackage[utf8]{inputenc}
\usepackage[english]{babel}

\usepackage{amsmath,amssymb,amsfonts,amsthm}
\usepackage[font=scriptsize,labelfont=bf]{caption}
\usepackage{mathrsfs}
\usepackage{textcomp}
\usepackage{graphics}
\usepackage{xcolor}
\usepackage[all,2cell]{xy}
\usepackage[ruled,linesnumbered, inoutnumbered]{algorithm2e}
\usepackage{tikz}
\usetikzlibrary{backgrounds,patterns,arrows,positioning,fit,shapes.geometric,decorations.pathmorphing}
\usepackage{url}
\usepackage[capitalize]{cleveref}
\Crefname{equation}{}{}

\newtheorem{theorem}{Theorem}
\newtheorem{lemma}{Lemma}
\newtheorem{proposition}[theorem]{Proposition}

\theoremstyle{definition}

\newtheorem{assumption}{Assumption}
\newtheorem{example}{Example}
\newtheorem{remark}{Remark}

\DeclareMathAlphabet{\mathpzc}{OT1}{pzc}{m}{it}

\newcommand{\Var}[1]{\mathcal{#1}}

\newcommand{\vect}[1]{\mathbf{#1}}

\newcommand{\Tang}[2]{\mathrm{T}_{#1} {#2}}
\newcommand{\Norm}[2]{\mathrm{N}_{#1} {#2}}

\newcommand{\diag}{\operatorname{diag}}

\newcommand{\vecc}[1]{\operatorname{vec}(#1)}
\newcommand{\R}{\mathbb{R}}
\newcommand{\deriv}[2]{\mathrm{d}_{#2}#1}

\DeclareMathOperator*{\argmin}{argmin}

\newcommand{\SFF}{\mathit{I\!I}}

\crefname{equation}{}{}
\crefname{equation}{}{}

\title[Sensitivity of low-rank matrix recovery]{Sensitivity of low-rank matrix recovery}

\author{Paul Breiding}
\thanks{PB: MPI MiS Leipzig,
Inselstr.\ 22, 04103 Leipzig, Germany. \texttt{paul.breiding@mis.mpg.de}.
PB has received funding from the European Research Council
(ERC) under the European Union's Horizon 2020 research and innovation programme (grant agreement No 787840) and from the Deutsche Forschungsgemeinschaft (DFG) -- Projektnummer 445466444.
}

\author{Nick Vannieuwenhoven}
\thanks{NV: KU Leuven, Department of Computer Science, Celestijnenlaan 200A, B-3001 Leuven, Belgium; Leuven.AI, KU Leuven Institute for AI, B-3000 Leuven, Belgium, \texttt{nick.vannieuwenhoven@kuleuven.be}. Supported by a Postdoctoral Fellowship of the Research Foundation---Flanders (FWO) with project 12E8119N}


\begin{document}

\begin{abstract}
We characterize the first-order sensitivity of approximately recovering a low-rank matrix from linear measurements, a standard problem in compressed sensing. A special case covered by our analysis is approximating an incomplete matrix by a low-rank matrix. This is one customary approach to build recommender systems.
We give an algorithm for computing the associated condition number and demonstrate experimentally how the number of linear measurements affects it.

In addition, we study the condition number of the rank-$r$ matrix approximation problem. It measures in the Frobenius norm by how much an infinitesimal perturbation to an arbitrary input matrix is amplified in the movement of its best rank-$r$ approximation. We give an explicit formula for the condition number, which shows that it does depend on the relative singular value gap between the $r$th and $(r+1)$th singular values of the input matrix.
\end{abstract}

\maketitle

\section{Introduction} \label{sec_introduction}

\textit{Compressed sensing} \cite{CRT2006, Donoho2006,DE2011,EK2012,FR2013} is a general methodology for recovering an unknown but structured signal $y \in \R^k$ from a measurement $a = L(y) \in \R^\ell$, where $\ell$ can be much smaller than $k$ and $L$ is a sensing operator.
The goal is to recover the unknown signal using only information about the compressed signal.
We consider only affine linear maps as sensing operators in this paper.

\textit{Low-rank matrix recovery} is a specific instance of compressed sensing. Herein, it is assumed that the unknown signal, an $m \times n$ matrix $Y$, (approximately) exhibits a low-rank structure of known rank $r$. The goal is to find a rank-$r$ matrix close to the unknown matrix $Y$ from the compressed sensing $A=L(Y)$.

A prominent application of low-rank matrix recovery is in \textit{collaborative filtering} and \textit{recommender systems}. Consider the so-called \emph{Netflix problem}~\cite{BL2009} for instance. Here, the data consists of an $m \times n$ matrix for $m$ users and $n$ movies and the $(i,j)$th entry contains the rating of user $i$ for movie $j$. Not all users have rated every movie. Thus, not all entries of the data matrix are available; it is incomplete. Filling in the missing values corresponds to predicting personalized movie ratings for each user. A common assumption is that the rating of movies by users is determined by unobserved latent factors, and that a low-rank factorization reveals these factors. This assumption was exploited by several submissions of the Netflix prize competition, including SVD++ \cite{Koren2008}, timeSVD++ \cite{Koren2009b}, and the eventual winning solution \cite{Koren2009}.
Recovering these latent factors from incomplete observations is a low-rank matrix recovery problem. Indeed, if the number of known ratings is $\ell$, we can arrange the entries of this incomplete matrix in a vector $A\in\mathbb R^\ell$. The projection from matrices to incomplete matrices is then a (linear) coordinate projection $L:\mathbb R^{m\times n}\to\mathbb R^\ell$.

Problems that can be solved with low-rank matrix recovery include collaborative filtering \cite{BKV2009, RS2005}, image inpainting  \cite{Magazine2014,KSK2015,MatrixInpainting2015}, dimensionality reduction \cite{SW2006, SY2007}, embedding problems \cite{LLE2995}, and multi-class learning \cite{AEP2008, OTM2010}. In all of these applications it is important to understand the sensitivity of the output with respect to perturbations in the input.
Eisenberg \cite{Eisenberg2020} summarizes this as follows:
\begin{quote}
  \emph{```Many investigations of big data solve inverse problems [...].
  The sensitivity of results to uncertainties [...] is crucial to determine the reliability and thus utility of results.''}
\end{quote}
For the Netflix problem this translates to the pertinent question of how sensitive the predicted ratings are to small perturbations in the known ratings (which by their nature are never truly exact).
The urgency of this question is underlined by recent concerns about reproducibility \cite{FBCJ2021}, well-posedness \cite{TTV2019}, and sensitivity \cite{AD2021} in recommender system technologies. For example, \cite{FBCJ2021} reported that the results of less than half of the considered conference papers could be reproduced. A priori, one potential, source of divergence in predicted ratings can be due to slight differences in the way the input data is centered, e.g., rounding to double, single, or half precision floating-point numbers. This may seem insignificant, but it is known that approximation problems can be very sensitive to changes to the input data due to ill conditioning \cite{BV2020}.
How can we ascertain whether these small perturbations are not propagated to large proportions in the final predicted ratings in recommender systems based on low-rank matrix factorization?

In this paper, we consider foregoing question for low-rank matrix recovery in the setting where $L$ can be any affine linear map. Formally, the low-rank matrix recovery problem consists of solving the nonlinear least-squares problem
\begin{equation}\label{lrr_problem}\tag{R}
\argmin_{\substack{Y\in \mathbb R^{m\times n},\; \mathrm{rank}(Y)=r}}\ \;\frac{1}{2} \Vert A-L(Y)\Vert^2,
\end{equation}
where $\Vert \cdot \Vert$ is the Euclidean norm on $\mathbb R^\ell$. We assume that we are given a \textit{well posed} problem instance. That is, a solution exists, is unique, and is locally continuous. We will return to discuss this assumption in \cref{sec:contributions}. For now, it suffices to know that if $\ell > r(m+n-r)$, then for almost all (affine) linear maps $L$ and almost all incomplete matrices $A$, the least-squares problem \cref{lrr_problem} is well posed by \cite[Q\&A 7]{BGMV2021}.

\subsection*{Contributions}

To answer the previous question for low-rank matrix recovery, we characterize the (first-order) sensitivity of the output, the unknown low-rank matrix $Y\in\mathbb R^{m\times n}$, with respect to small perturbations of the input data, the compressed sensing $A\in\mathbb R^\ell$. For this, we compute the \emph{condition number} $\kappa_\mathrm{recovery}(A,Y)$ of the nonlinear least-squares problem \cref{lrr_problem}. We give a formal definition of this number in \cref{sec:contributions}, but at this point it suffices to think about an asymptotically sharp bound
\[
\Vert Y - Y'\Vert_F \leq \kappa_\mathrm{recovery}(A,Y)\,\Vert A-A'\Vert,
\]
where $\| \cdot \|_F$ denotes the Frobenius norm.
Here, $Y'$ is the solution of \cref{lrr_problem} for the input $A'=A+\Delta A$, which is a small perturbation of $A$. By asymptotically sharp we mean that the inequality is a sharp inequality in the limit as $\Vert \Delta A\Vert \to 0$. We stress that the foregoing bound holds irrespective of the specific algorithm that is employed to obtain the low-rank matrix $Y$. It is an intrinsic property of the low-rank matrix recovery problem, a measure of its numerical hardness \cite{BCSS,BC2013}.

Our first main contribution is a numerical linear algebra algorithm for computing the condition number $\kappa_\mathrm{recovery}(A,Y)$ of low-rank matrix recovery. This algorithm is presented in \cref{alg_matrix_recovery}. We show in \cref{prop_complexity} below that for certain structured sensing operators, including coordinate projections, the computational complexity of the algorithm is $\mathcal{O}( \phi s^3 )$, where $s = (m+n-r)r$ is the problem size\footnote{The problem size is defined here as the dimension of the optimization domain in \cref{lrr_problem}.
It will be stated formally in \cref{sec:contributions}.} and $\ell = \phi s$, where the \emph{oversampling rate} $\phi>1$ is typically a small constant (up to a factor $\phi$, this is the same complexity as one step of a standard (Riemannian) Newton method \cite{AMS2008,Boumal2020} for solving optimization problem \cref{lrr_problem}).  We apply the algorithm in \cref{sec:experiments} to small-scale matrix recovery problems to assess the impact of oversampling ($\phi>1$) on the condition number.

Our second contribution is an explicit formula of the condition number for the special case of \emph{low-rank approximation}. Here, the input data is $A\in\mathbb R^{m\times n}$ and the problem is approximating $A$ with a matrix of low rank $r$, i.e., solving
\begin{align}\label{lra_problem} \tag{A}
\argmin_{\substack{Y\in \mathbb R^{m\times n},\; \mathrm{rank}(Y)=r}}\ \;\frac{1}{2}\Vert A-Y\Vert^2_F,
\end{align}
where the norm is the Frobenius norm.
This problem is the special case of \cref{lrr_problem} when $L$ is the identity map. The usual approach for solving the low-rank approximation problem is by computing a compact singular value decomposition (SVD) $A = \sum_{i=1}^{\min\{m,n\}} \sigma_i \vect{u}_i \vect{v}_i^T$ with $\sigma_1\geq \cdots\geq \sigma_{\min\{m,n\}} \geq0$.
Then, a solution of \cref{lra_problem} is given by the truncated SVD $Y=\sum_{i=1}^{r} \sigma_i \vect{u}_i \vect{v}_i^T$.
We denote the condition number in this case by $\kappa_\mathrm{approximation}(A,Y)$. Our second main result, \cref{main_thm1}, characterizes the condition number of low-rank approximation of $A$:
\begin{align} \label{eqn_main_thm1}
\kappa_\mathrm{approximation}(A,Y) = \frac{1}{1-\frac{\sigma_{r+1}}{\sigma_r}} = \frac{\sigma_r}{\sigma_r - \sigma_{r+1}}.
\end{align}
The sensitivity of approximating $A$ with the low-matrix~$Y$ thus depends on the \emph{singular value gap} between $\sigma_{r+1}$ and $\sigma_r$ of $A$. The input $A$ is \emph{ill-posed}, if $\sigma_{r-1} = \sigma_r$.

The last statement might seem contradictory to some literature, like \cite{DI2019}, that might be interpreted as suggesting that ``low-rank matrix approximations do not need a singular value gap'' to have a small condition number. An informal example makes it clear, however, that \textit{if the recovered low-rank matrix is of interest, rather than the approximation error}, then a singular value gap is required for a meaningful interpretation. For the example we denote by $e_1, e_2\in\mathbb{R}^2$ the two standard basis vectors and take $0 \le \epsilon < 1$. The best rank-$1$ approximation of $A=\left[\begin{smallmatrix}1+\epsilon&0\\0&1-\epsilon\end{smallmatrix}\right]$ is
\[
 Y = (1+\epsilon) e_1 e_1^T = \begin{bmatrix} 1 &0 \\ 0& 0 \end{bmatrix} + \epsilon \begin{bmatrix} 1 & 0\\ 0& 0 \end{bmatrix},
\]
and we have $\kappa_{\mathrm{approximation}}(A,Y) = \frac{1}{2}(1 + \epsilon^{-1})$.
Therefore, a large deviation of $Y$ may be expected when perturbing $A$. For example, perturbing $A$ by $\epsilon(e_1 e_2^T + e_2 e_1^T)$, results in the matrix~$A' := \left[\begin{smallmatrix} 1 + \epsilon & \epsilon \\ \epsilon & 1 - \epsilon \end{smallmatrix}\right]$.
The best rank-$1$ approximation of $A'$ is
{\small\begin{equation*}
Y'
= \frac{1+\sqrt{2}\epsilon}{2(2 + \sqrt{2})} \begin{bmatrix} 1 + \sqrt{2} \\ 1 \end{bmatrix} \begin{bmatrix} 1 + \sqrt{2} \\ 1 \end{bmatrix}^T
= \begin{bmatrix}
\frac{1}{2} + \frac{\sqrt{2}}{4} & \frac{1}{2\sqrt{2}} \\[5pt]
\frac{1}{2\sqrt{2}} & \frac{1}{2(2+\sqrt{2})}
\end{bmatrix} + \epsilon \begin{bmatrix}
\frac{1}{2} + \frac{1}{\sqrt{2}} & \frac{1}{2} \\[5pt]
\frac{1}{2} & \frac{1}{4(1+\sqrt{2})}
\end{bmatrix}.
\end{equation*}}%
Hence, a unit-order change results between $Y'$ and $Y$ from a perturbation of size $\sqrt{2}\epsilon$, as could have been anticipated from $\kappa_\mathrm{approximation}(A,Y) \approx \epsilon^{-1}$. Another more formal example is given in Example \ref{Example} in \cref{sec_low_rank_approximation} below.

\subsection*{Outline}
The outline of this paper is as follows. In the next section, we compare and contrast the perhaps surprising result for low-rank approximation to existing insights from the literature. Thereafter, \cref{sec:contributions} formally states the main results and assumptions of our study.  \Cref{sec:H} investigates the Hessian of the objective function from \cref{lra_problem}, which provides a crucial contribution to the condition numbers of both problems \cref{lrr_problem,lra_problem}. Armed with insights about the Hessian, we characterize the condition number of \cref{lra_problem} in \cref{sec_low_rank_approximation}. The condition number of \cref{lrr_problem} is analyzed in \cref{sec:recovery}; in this case, we are unfortunately not able to derive a closed expression. For this reason, \cref{alg_matrix_recovery} presents a numerical algorithm for computing it. Numerical experiments with both low-rank approximation and recovery are featured in \cref{sec:experiments}.

\subsection*{Acknowledgements} We thank Sebastian Kr\"amer for valuable feedback that led to several improvements in the presentation of our results. In particular, formulating \cref{Example} by formalizing the introductory example was suggested by him.

\section{Comparison to prior results}\label{sec:Hackbusch}

In the literature we did not find results on the sensitivity of low-rank matrix recovery.
However, for the special case of low-rank approximation there are several. Some of them might seem to contradict our result, while the final one corroborates it. For this reason, we carefully discuss the prior literature.

\subsection{Drineas and Ipsen's no-gap result}
Drineas and Ipsen's article \cite{DI2019} is titled \textit{``Low-rank matrix approximation do not need a singular value gap.''} At first sight, this seems to contradict our results, but on closer inspection the paradox quickly disappears. Dirineas and Ipsen study error bounds for the \emph{approximation error} of a low-rank approximation, as measured by the Schatten $p$-norm of the residual $P_U^\perp A$, where $A \in \R^{m \times n}$ is the matrix to approximate and $P_U^\perp$ projects onto the orthogonal complement of a fixed $r$-dimensional subspace $U \subset \R^m$. That is, they derive error bounds for $\| P_U^\perp A \|_p$ as either the fixed subspace $U$ or the matrix $A$ is perturbed in \cite[Theorem 1]{DI2019} and \cite[Theorem 2]{DI2019}, respectively. For example, when perturbing $A$, \cite[Theorem 2]{DI2019} states that
\[
 \big| \| P_U^\perp A \|_p - \| P_U^\perp A' \|_p \big| \le \|A - A'\|_p.
\]
Our results, on the other hand, describe what happens to the best rank-$r$ approximation of $A$ as it is perturbed. That is, using terminology closer to \cite{DI2019}, we show that
\[
\| P_{U^*}^\perp(A) \, A - P_{U^*}^\perp(A') \, A' \|_F \le \frac{\sigma_r}{\sigma_r - \sigma_{r+1}} \| A - A' \|_F ,
\]
where $P_{U^*}^\perp(A)$ projects $A$ to the best rank-$r$ approximation.\footnote{Equivalently, but closer to \cite{DI2019} in formulation, it projects the column space of $A$ to $U^*$, the $r$-dimensional subspace of left singular vectors associated to the largest $r$ singular values.}
Note that in our result both $A$ and the projector $P_{U^*}^\perp(A)$ are perturbed  as $P_{U^*}^\perp(A)$ varies with $A$.

\subsection{An error bound of Hackbusch}
Next, we discuss the result from \cite{Hackbusch2016} by Hackbusch.
For this we let $A\in\mathbb R^{m\times n}$ and we denote by $A':=A+\Delta A$ a perturbation of~$A$. If $Y'$ is the best rank-$r$ approximation of $A'$ and if $Y_\mathrm{computed}$ is any other rank-$r$ matrix, then Theorem~4.5 in \cite{Hackbusch2016} asserts that
\begin{equation}\label{hackbusch}
  \Vert Y'-Y_\mathrm{computed}\Vert_F \leq q\,\Vert A'-Y_\mathrm{computed}\Vert_F,
\end{equation}
where $q=\tfrac{1+\sqrt{5}}{2}\approx 1.62$ is a constant that does not depend on the input data.

The main rationale for this bound is that $Y_\mathrm{computed}$ could be a cheap approximation of the rank-$r$ truncated SVD, e.g., obtained from randomized methods \cite{HMT2011} or adaptive cross approximation \cite{BR2000}.
The bound states that if the approximation is $Y_\mathrm{computed}$, then $Y_\mathrm{computed}$ deviates from $Y'$ by at most $q$ times $\Vert A'-Y_\mathrm{computed}\Vert_F$. The latter can be computed from the data, so the quality of the computation can be assessed.

The fact that \cref{hackbusch} involves a constant upper bound seems contradictory to \cref{eqn_main_thm1}. However, our result states that for sufficiently small $\Vert A' - A\Vert_F$ we have
\begin{equation}\label{hackbusch2}
  \Vert Y' - Y \Vert_F \leq \frac{\sigma_r}{\sigma_r- \sigma_{r+1}}\,\Vert A' - A \Vert_F + o( \Vert A' - A \Vert_F^2 ),
\end{equation}
where $Y$ is a best rank-$r$ approximation of $A$, and $\sigma_r,\sigma_{r+1}$ are the $r$th and $(r+1)$th singular values of $A$.
This means that a small perturbation $\Delta A$ of the input $A$  is amplified in the output by the condition number $\kappa_\text{approximation}(A,Y)$ in the worst case. If $0 \ne \sigma_{r}\approx\sigma_{r+1}$, this factor is huge.

Assume that $A$ is the true matrix we want to compute a rank-$r$ approximation of and that $A'=A+\Delta A$ is a perturbation of $A$, e.g., due to roundoff or measurement errors. In this case, the bound \cref{hackbusch} does not tell the whole story and could be complemented with \cref{eqn_main_thm1}. Indeed, even if we can approximate $A'$ closely by $Y_\text{computed}$ so that $\Vert A'-Y_\mathrm{computed}\Vert_F$ is small, the matrix $Y_\mathrm{computed}$ can still be far from the best rank-$r$ approximation $Y$ of the true matrix $A$. Combining Hackbusch's result with \cref{eqn_main_thm1} yields
\begin{align*}
\Vert Y-Y_\mathrm{computed}\Vert_F
&\leq \Vert Y-Y'\Vert_F + \Vert Y'-Y_\mathrm{computed}\Vert_F \\
&\leq q \Vert A'-Y_\mathrm{computed}\Vert_F + \frac{\sigma_r}{\sigma_r-\sigma_{r+1}} \Vert A' - A \Vert_F.
\end{align*}
The first term of the final bound follows from Euclidean geometry, while the second term is the effect of curvature of the manifold of rank-$r$ matrices.

Finally, observe that both \cref{hackbusch} and \cref{hackbusch2} agree on a constant upper bound when~$A'$ is a perturbation of the rank-$r$ matrix $A = Y_\text{computed}$. In this case, $\sigma_r > \sigma_{r+1} = 0$ so that $\kappa_\text{approximation}(A,Y_\text{computed}) = 1$.

\subsection{First-order perturbations of the SVD by Hua and Sarkar}
The earliest result on the sensitivity of the best low-rank approximation to a matrix we could locate in the literature is by Hua and Sarkar \cite{HS89}. They show that ``the first-order perturbations in the SVD truncated matrices [...] can be simply expressed in terms of the perturbations in the original data matrices'' and they conclude from their analysis that ``the SVD truncations do not affect the first order perturbations''. This also seems to contradict \cref{eqn_main_thm1}, where we show that a best rank-$r$ approximation \emph{can} change by (much) more than the norm of the perturbation.

The paradox disappears when we take into account that Hua and Sarkar assume that the input $A$ is itself a rank-$r$ matrix. Thus, the $(r+1)$th singular value of $A$ is $\sigma_{r+1}=0$ and so, by \cref{eqn_main_thm1}, we have, once more, $\kappa_\mathrm{approximation}(A,Y) = 1$, which is fully consistent with their result.

\subsection{Perturbation expansions of Vu, Chunikhina, and Raich}
The main result of Vu, Chunikhina, and Raich \cite[Theorem 1]{VCR2021} turns Feppon and Lermusiaux's analysis \cite{FL2018} into a rigorous perturbation bound for the best rank-$r$ approximation for arbitrary input matrices. We also use Feppon and  Lermusiaux's work in \cref{sec_low_rank_approximation,sec:recovery}.
Consequently, the effect of curvature pops up in \cite[Theorem 1]{VCR2021}, consistent with \cref{eqn_main_thm1}. Nevertheless, we think our first-order error bound
\[
 \| Y - Y' \|_F \le \frac{\sigma_r}{\sigma_r - \sigma_{r+1}} \| A - A' \|_F + \mathcal{O}(\|A-A'\|_F^2),
\]
where $A'$ is a perturbation of $A$ and $Y$ and $Y'$ are the best rank-$r$ approximations of $A$ and $A'$ respectively, is more succinct than the bound in \cite[Theorem~1]{VCR2021}. In addition, our analysis extends to the low-rank matrix recovery problem.

\section{Statement of the main results}\label{sec:contributions}
As in the introduction we consider an affine linear map
\begin{align}\label{eqn_sensing_op}\tag{S}
L:\mathbb R^{m\times n}\to \mathbb R^\ell, \, Y\mapsto M(Y) + b,
\end{align}
which is called the sensing operator.
Let
$$\Var{M}_r:=\{X\in\mathbb R^{m\times n} \mid \mathrm{rank}(X)=r\}$$ be the set of matrices of rank equal to $r$. It is a \textit{smooth embedded submanifold} of $\R^{m\times n}$ of dimension $\dim \Var{M}_r = (m+n-r)r$ \cite{HM1994}. This implies that the set $\Var{M}_r$ is equipped with a topology and smoothness structure that is inherited from the ambient space~$\mathbb{R}^{m\times n}$. This enables a vast generalization of calculus on such domains \cite{Lee2013}. Precisely this smooth structure will make it much easier to compute the desired condition numbers. The set of sensed matrices will be denoted by
\[
\Var{I}_r:=L(\Var{M}_r).
\]
We also define the set of matrices of rank bounded by $r$:
\[
\Var{M}_{\leq r}:=\{X\in\mathbb R^{m\times n} \mid \mathrm{rank}(X)\leq r\}.
\]
It is both the Euclidean closure of $\Var{M}_r$ and a real algebraic variety in $\R^{m\times n}$, defined by the vanishing of $(r+1)\times (r+1)$-minors \cite{Harris1992}.

The goal of this paper is to determine the first-order sensitivity of \cref{lrr_problem}. The input to this problem is a compressed sensing $A \in \R^\ell$, while the output is necessarily restricted to be a rank-$r$ matrix.
The first complication one encounters is that there can be no or several solutions $Y$ for an input $A$. A priori we should expect to deal with a
\emph{set-valued solution map}
$
R : \R^{\ell} \rightrightarrows \Var{M}_r, \; A \mapsto
\argmin_{Y \in \Var{M}_r} \, \frac{1}{2} \Vert A-L(Y)\Vert^2.
$
The condition number of this map can be analyzed with the general techniques we introduced in \cite{BV2020}. In low-rank recovery, however, the geometry of the problem is more well-behaved than the general case. This allows for a clearer presentation that eliminates the intricacy of solution manifolds in \cite{BV2020}. We explain this next.

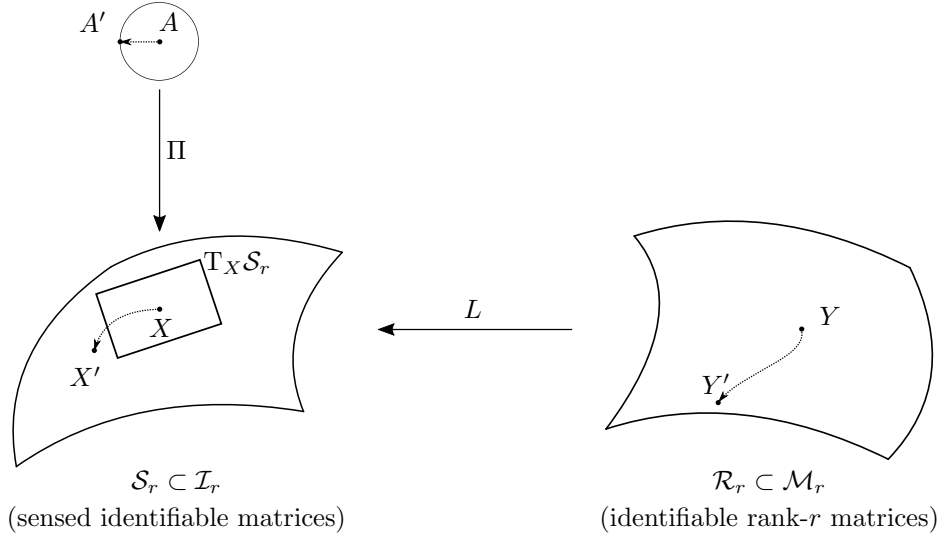
\begin{figure}[tb]
\begin{tikzpicture}[remember picture]
\node[inner sep=0] at (0,0)
{\includegraphics[width=.98\textwidth]{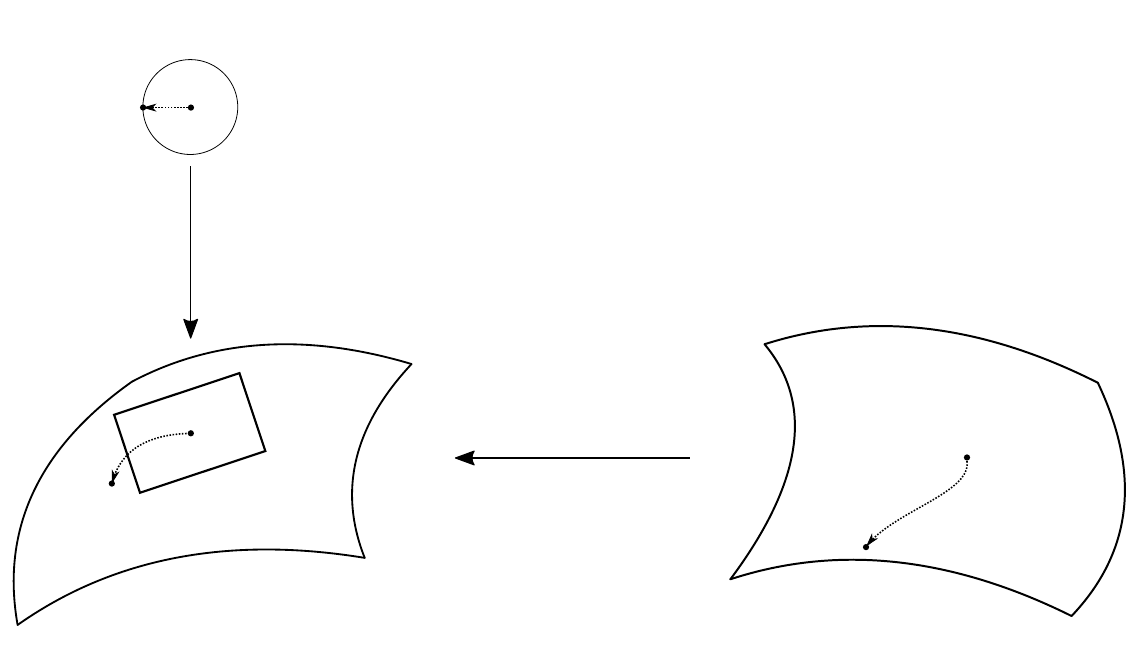}};
\node at (-4,2.6) {$A$};
\node at (-5.0,2.6) {$A'$};
\node at (-3.9,0.9) {$\Pi$};
\node at (-4.1,-1.45) {$X$};
\node at (-3.1,-0.6) {$\Tang{X}{\Var{S}_r}$};
\node at (-5.1,-2.1) {$X'$};
\node at (4.7,-1.3) {$Y$};
\node at (3.2,-2.2) {$Y'$};
\node at (0,-1.2) {$L$};
\node at (-3.9,-3.5) {$\mathcal S_r\subset \mathcal I_r$};
\node at (-3.9,-4) {(sensed identifiable matrices)};
\node at (3.9,-3.5) {$\mathcal R_r\subset \mathcal M_r$};
\node at (3.9,-4) {(identifiable rank-$r$ matrices)};
\end{tikzpicture}
\caption{\label{fig_geometry} The simplified geometry of this paper. On the right is the submanifold $\mathcal R_r\subset \mathcal M_r$ of identifiable rank-$r$ matrices, and on the left is the submanifold $\mathcal S_r\subset \mathcal I_r$ of sensed identifiable matrices. If the affine linear map $L$ is generic, then it restricts to a diffeomorphism $\mathcal R_r\to\mathcal S_r$. The low-rank matrix recovery problem consists of two steps: (i) projecting the data point $A\in\R^\ell$ to $X\in\Var{S}_r$ with $\Pi$, and (ii) finding $Y\in\Var{R}_r$ with $L(Y)=X$. Therefore, the sensitivity of the output $Y$ with respect to the input perturbation $A' - A$ depends on the combined impact of (i) the curvature of $\Var{S}_r$ which causes $X$ to move to $X' \in \Var{S}_r$ as $A$ moves to $A'$, and (ii) the sensitivity of inverting $L\mid_{\mathcal R_r}$ which forces $Y$ to move to $Y' \in \Var{R}_r$ as $X$ moves to $X'$.}
\end{figure}

The geometry of our setting is depicted in \cref{fig_geometry}. It illustrates that \cref{lrr_problem} decouples into two subproblems:
\begin{enumerate}
\item[(i)] minimizing the distance from $\mathcal I_r$ to $A$, and
\item[(ii)] inverting the map $L$.
\end{enumerate}
Fortunately, under a mild assumption, $L^{-1}$ is a differentiable \textit{function} almost everywhere in the precise sense of \cref{prop_ass1} below. This assumption is the following.
\begin{assumption}\label{ass1}
We assume that $L(Y)=M(Y)+b$ can be \emph{generically identified}.
\end{assumption}
Being generically identifiable means that there is a Zariski open algebraic subvariety $\Sigma\subset \Var M_{\leq r}$, such that~$L^{-1}(L(Y)) \cap \Var M_r = \{Y\}$ for all $Y\in\Var M_r\setminus \Sigma$. Q\&A 7 in \cite{BGMV2021} shows that almost all $L$ have this property if $\ell \ge \dim \Var{M}_r = (m+n-r)r$.

The next result is \cite[Q\&A 12]{BGMV2021}.
\begin{proposition}\label{prop_ass1}
Under \cref{ass1} there exist smooth embedded submanifolds $\Var{R}_r \subset \Var{M}_r$ and $\Var{S}_r \subset \Var{I}_r = L(\Var{M}_r)$, which are both dense in their supsets, such that
\[
L|_{\Var{R}_r} : \Var{R}_r \to \Var{S}_r
\]
is a global diffeomorphism.
\end{proposition}
A global diffeomorphism is a smooth bijective map between manifolds whose inverse map is smooth.

%

For minimizing the distance from $\Var{S}_r$ to the input $A\in\R^\ell$ we consider the following open subset of the input space $\R^\ell$:
\[
\Var{D} = \R^\ell \setminus \Var{Z},
\]
where $\Var{Z}$ is (the Euclidean closure of) the set of points for which $\min_{X\in\Var{S}_r}\,\tfrac{1}{2}\Vert A-X\Vert^2$ does not have a \textit{unique} solution (because it has multiple or no solutions). Note that we have replaced $\Var{I}_r$ by the submanifold $\Var{S}_r$ here. Since $\Var{S}_r$ is an embedded\footnote{This means that its smooth structure is compatible with the smooth structure on $\R^\ell$.%
} submanifold of $\R^\ell$, the existence of a \textit{tubular neighborhood} \cite{Lee2013} of $\Var{S}_r$, an open neighborhood containing $\Var{S}_r$ in $\R^\ell$, guarantees that $\mathcal D$ contains at least this open subset.\footnote{Erd\"os's result~\cite{Erdos1945} implies that the set of points which have several infima of the distance function to $\mathcal{S}_r$ is of Lebesgue measure zero.
However, $\mathcal{S}_r$ may not be closed, and so there can be points with no minimizer on $\mathcal{S}_r$.
The set of such points can even be full-dimensional. Think of a cusp with the node removed. Recently, explicit examples of nonclosedness were presented in \cite{TTV2019} in the context of matrix completion.}
On $\Var{D}$ we can define $\Pi: \Var{D} \to \Var{S}_r, A\mapsto \argmin_{X\in\mathcal S_r}\,\tfrac{1}{2}\Vert A-X\Vert^2$, the projection onto $\mathcal S_r$.

In summary, the foregoing closer look at the geometry of problem \cref{lrr_problem} allows us to arrive at a \emph{recovery map} $R = (L|_{\mathcal R_r})^{-1}\circ \Pi$. This is the map
\begin{equation} \label{eqn_vandereycken_formulation2}
 R: \mathcal D \to \Var{R}_r,\quad A \mapsto \argmin_{Y \in \Var{R}_r} \frac{1}{2} \| A - L(Y) \|^2.
\end{equation}
It is a smooth (uni-valued) map!

With the foregoing concessions (\cref{ass1}, open dense submanifolds $\Var{R}_r$ and~$\Var{S}_r$, removing $\mathcal Z$ from the domain) we can apply Rice's \cite{Rice1966} classic definition of the condition number of a map for $A \in \mathcal D$:
\begin{equation}\label{def_cond}
\kappa_\mathrm{recovery}(A, Y) = \lim_{\epsilon \to 0}\sup_{\substack{\Delta A\in\mathbb R^\ell,\\ \Vert \Delta A\Vert \leq \epsilon}}\,\frac{\Vert R(A) - R(A+\Delta A)\Vert_F}{\Vert \Delta A\Vert},
\end{equation}
where $Y = R(A) \in \Var{R}_r$ is the recovered rank-$r$ matrix. If $A\in\mathcal Z$ is outside the locus, where we can obtain the recovery map (\ref{def_cond}), we define $\kappa_\mathrm{recovery}(A, Y) = \infty$.
\begin{remark}\label{rmk1}
Note that \cref{eqn_vandereycken_formulation2,def_cond} allow us to study the condition number of the global minimizer $Y = R(A)$ of \cref{eqn_vandereycken_formulation2}. By considering the graph of $L : \Var{R}_r \to \Var{S}_r$, we see that the results of \cite{BV2020} also apply in their general form. This means that the analysis in this paper also covers local minima and critical points $Y$ with $A - L(Y) \perp \Tang{L(Y)}{\Var{S}_r}$ as in \cite{BV2020}. Nevertheless, for concreteness we focus on (local) minima, because they are the main interest in applications.
\end{remark}

\subsection{Low-rank matrix recovery}
The condition number of (well-posed) low-rank matrix recovery in \cref{eqn_vandereycken_formulation2} at $A\in \mathcal D$ with output $Y=R(A)$ can be obtained from Theorem 7.3 in \cite{BV2020}:
\begin{equation}\label{def_CN}\tag{C}
\kappa_\mathrm{recovery}(A,Y) = \| (M|_{\mathrm T_Y\mathcal M_r})^{-1} H_{A,X}^{-1} \|_2.
\end{equation}
Herein, $\| \cdot \|_2$ is the spectral norm relative to the Frobenius norms on $\R^{m\times n}$ and the Euclidean norm $\R^\ell$, $M|_{\mathrm T_Y\mathcal M_r}$ is the derivative of $L|_{\Var{M}_r}$, and $H_{A,X}$ is the \emph{Riemannian Hessian} of the squared distance function $d_A : \Var{S}_r \to \R, X \mapsto \frac{1}{2} \| A - X \|^2$ at the point~$X=L(Y)$.
This Riemannian Hessian generalizes the classic Euclidean Hessian and contains the second derivatives of $d_A$ on $\Var{S}_r$; it is discussed in \cref{sec:H}.

\begin{figure}
    \begin{tikzpicture}[scale = 1.54]
      \draw[black]   plot[smooth,domain=-1.2:1.2] (\x, {-(\x)*(\x)});
      \draw[black]  (0,-1/2) circle[radius=1/2];

      \draw[black]  (0,0) circle[radius=1.25pt] node[above, xshift=5pt] {$X$};
      \draw[black]  (-0.35,-0.35 * 0.35) circle[radius=1.25pt] node[above, xshift=-15pt] {$X+\Delta X$};

      \draw[->] (0,-0.3)  -- (0,-0.1);

      \fill[gray]  (0,-1/2) circle[radius=1.25pt];
      \fill[black]  (0,-0.4) circle[radius=1.25pt] node[right] {$A$};

      \draw[->] (-0.1,-0.4)  -- (-0.19,-0.4)  node[left, fill=white] {$A+\Delta A$};
      \draw[->, shorten >=3pt] (-0.2,-0.35)  -- (-0.35,-0.35 * 0.35) ;

           \fill[white]  (0,0) circle[radius=1.24pt];
      \fill[white]  (-0.35,-0.35 * 0.35) circle[radius=1.24pt];
 \end{tikzpicture}
\caption{\label{fig_osculating} The picture shows how curvature affects the sensitivity of computing closest points on nonlinear objects. In this case, the curvature of the parabola amplifies the error $\Delta A$ in $A$. The amplification of errors is determined by the eigenvalues of the Riemannian Hessian $H_{A,X}$.
}
\end{figure}
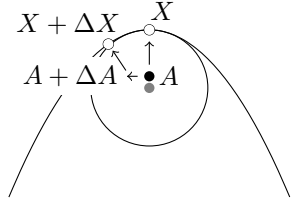

As one can see from \cref{fig_geometry} and also from the formula \cref{def_CN}, the condition number $\kappa_\mathrm{recovery}(A,Y)$ is determined by two parts:
\begin{enumerate}
\item[(i)] the sensitivity of the recovery map $(L|_{\mathcal M_r})^{-1}$, and
\item[(ii)] the \emph{curvature} of the manifold of sensed rank-$r$ matrices $\Var{S}_r$ at $L(Y)$.
\end{enumerate}
The effect of curvature on condition is depicted in \cref{fig_osculating}. If $A$ is a \emph{center of curvature} with base point $X$ of the parabola-shaped manifold, then the Riemannian Hessian $H_{A,X}$ is not invertible. In this case we have~$\kappa_\mathrm{recovery}(A,Y)=\infty$ and we call the input $A$ \emph{ill-posed}. The center of curvature for $X$ is shown in \cref{fig_osculating} as the gray point in the center of the displayed circle.

\subsection{Low-rank matrix approximation}\label{sec:contr_lra}
We turn to the special case when $L$ is the identity map. This corresponds to the problem of approximating a matrix $A \in \R^{m\times n}$ by a rank-$r$ matrix.

The condition number of \emph{low-rank approximation} is also given by \cref{def_CN}, where $L$ is the identity and $\Var{R}_r = \Var{S}_r$. Our main result in this setting is the following result.
\begin{theorem}\label{main_thm1}
Let $A = \sum_{i=1}^{\min\{m,n\}} \sigma_i \vect{u}_i \vect{v}_i^T$ be an SVD of $A$ with ordered singular values $\sigma_1\geq \ldots\geq \sigma_{\min\{m,n\}} \ge 0$. Let $Y\in R(A)$ be a rank-$r$ truncated SVD of $A$. Then, the condition number of finding a best rank-$r$ approximation at $(A,Y)$ is
\[
 \kappa_\mathrm{approximation}(A,Y) = \frac{1}{1 - \frac{\sigma_{r+1}}{\sigma_r}} = \frac{\sigma_r}{\sigma_r - \sigma_{r+1}}
\]
or $1$ if $\sigma_r = 0$.
\end{theorem}
We prove this theorem in \cref{sec_low_rank_approximation} below. Note that $\kappa_\mathrm{approximation}(A,Y)=\infty$ if and only if $\sigma_{r+1} = \sigma_r \ne 0$.


\section{The Riemannian Hessian of the distance function}\label{sec:H}

The \textit{Riemannian Hessian} \cite{Lee1997,riemannian_geometry,Petersen} generalizes the classic Hessian matrix from multivariate functions to maps on manifolds. We focus on the Riemannian Hessian of the distance function $d_A : \Var{X} \to \R, X \mapsto \frac{1}{2} \| X - A \|_F^2$ from $A \in \R^n$ to the smoothly embedded submanifold $\Var{X} \subset \R^n$.
This manifold is equipped with the \textit{Riemannian metric} inherited from the Euclidean space $\R^n$. That is, every tangent space $\Tang{X}{\Var{X}}$ is equipped with the inner product $g_X(x,y) = x^T y$, where $x, y \in \Tang{X}{\Var{X}}$ are viewed as vectors in $\R^n$.

The goal of this section is not to provide a rigorous derivation of the Riemannian Hessian in general, but rather present an accessible account for submanifolds of Euclidean space that highlights its connection to classic differential-geometric objects like the second fundamental form and Weingarten map, which will be used in the technical results. An alternative accessible account can be found in \cite[Chapter 5]{Boumal2020}.

Ignoring the manifold structure for a moment, in classic multivariate analysis the Hessian of $d_A$ at $X$ would be
\begin{align} \label{eqn_hessian} \tag{E}
 \deriv{(\deriv{d_A}{X})}{X}
 &= \deriv{(\dot{X} \mapsto \langle \dot{X}, X - A \rangle)}{X} \\
 \nonumber &= (\dot{X},\ddot{X}) \mapsto \langle \dot{X}, \ddot{X} \rangle - \langle (\deriv{\dot{X}}{X})(\ddot{X}), A - X \rangle.
\end{align}
In the classic setting, $\dot{X}$ is a vector in $\R^n$ which bears no particular relationship to~$X$. Hence, in Euclidean geometry, the second term involving $(\deriv{\dot{X}}{X})(\ddot{X})$ vanishes.

When $X$ is restricted to lie on a manifold $\Var{X}$, the interpretation of \cref{eqn_hessian} changes substantially. The derivative of a smooth map $f : \Var{X} \to \Var{Y}$ between manifolds at $X \in \Var{X}$ is a linear map $\deriv{f}{X} : \Tang{X}{\Var{X}} \to \Tang{f(X)}{\Var{Y}}$ between the respective tangent spaces \cite{Lee2013}. This means that $\dot{X}$ and $\ddot{X}$ are elements of $\Tang{X}{\Var{X}}$, which we can view in $\R^n$ as an affine linear space attached at $X$. Consequently, $(\deriv{\dot{X}}{X})(\ddot{X})$ should be interpreted as the directional derivative of the tangent vector $\dot{X}\in\Tang{X}{\Var{X}}$ as the base point $X \in \Var{X}$ is infinitesimally moved in the direction of $\ddot{X}\in\Tang{X}{\Var{X}}$. Based on this interpretation, circumventing vector fields, we could define
\[
 \nabla^2 : \Tang{X}{\Var{X}} \times \Tang{X}{\Var{X}} \to \R^n, \quad (\dot{X}, \ddot{X}) \mapsto \frac{\mathrm{d}}{\mathrm{d} t}\big|_{t=0} \mathrm{P}_{\Tang{\gamma(t)}{\Var{X}}} (\dot{X}),
\]
where $\gamma(t) \subset \Var{X}$ is an \textit{integral curve} \cite{Lee2013} realizing $\ddot{X}$, i.e., $\gamma$ is a smooth map from a neighborhood of $0 \in \R$ with $\gamma(0) = X$ and $\gamma'(0) = \ddot{X}$.

Since $\Tang{X}{\Var{X}}$ can be viewed as an affine linear subspace of $\R^n$, we can decompose the latter as $\R^n = \Tang{X}{\Var{X}} \oplus \mathrm{N}_{X}\Var{X}$, where $\mathrm{N}_{X}\Var{X}$ is the normal space of $\Var{X}$ at $X$, i.e., the orthogonal complement of $\Tang{X}{\Var{X}}$. Projecting $\nabla^2$ to the normal space yields a fundamental object in Riemannian geometry, the \textit{second fundamental form} $\SFF_X$ \cite{Lee1997,ONeill1983,ONeill2001,riemannian_geometry,Petersen}.
This is the bilinear map
\begin{equation}\label{SFF}
\SFF_X : \Tang{X}{\Var X}\times \Tang{X}{\Var{X}} \to \Norm{X}{\Var{X}},\quad (\dot{X}, \ddot{X}) \mapsto \mathrm{P}_{\Norm{X}{\Var{X}}}\left( \nabla^2(\dot{X}, \ddot{X}) \right).
\end{equation}
If we contract the output of this map with a normal vector $N \in \Norm{X}{\Var{X}}$, we obtain the so-called \textit{Weingarten map} or \textit{shape operator}, another classic and well-studied object in Riemannian geometry  \cite{Lee1997,ONeill1983,ONeill2001,riemannian_geometry,Petersen}:
\[
 S_N : \Tang{X}{\Var X} \times \Tang{X}{\Var X} \to \R, \quad (\dot{X}, \ddot{X}) \mapsto \langle N, \SFF_X(\dot{X}, \ddot{X}) \rangle.
\]
This Weingarten map is a self-adjoint operator $\Tang{X}{\Var{X}}\to \Tang{X}{\Var{X}}$.

At a critical point $N = A - X \in \Norm{X}{\Var{X}}$, it follows that we can interpret \cref{eqn_hessian} in terms of classic, well-studied objects in Riemannian geometry, namely as the following linear endomorphism on $\Tang{X}{\Var{X}}$:
\begin{equation}\tag{H}\label{eqn_Hess_distance}
 H_{A,X} = \mathbf{1}_{\Tang{X}{\Var{X}}} - S_{N},
\end{equation}
where $\mathbf{1}_{\Tang{X}{\Var{X}}}$ is the identity on ${\Tang{X}{\Var{X}}}$ and $S_N$ is viewed as linear map $\Tang{X}{\Var{X}} \to \Tang{X}{\Var{X}}$. The map $H_{A,X}$, or its matrix representation, is called the Riemannian Hessian.\footnote{For arbitrary $X$, the Riemannian Hessian is also defined by \cref{eqn_Hess_distance} for $N = \mathrm{P}_{\Norm{X}{\Var{X}}}(A - X)$ \cite{Petersen}.}

\subsection{Principal curvatures}
The Riemannian Hessian of $d_A$ contains geometric information about the way $\Var{X}$ curves inside of $\R^n$ \cite{Petersen}.
Let $N \in\mathrm{N}_X\Var X$, $\eta = \frac{N}{\Vert N\Vert}$, and $s:=\dim \Var{X}$.
The real eigenvalues $\lambda_1, \ldots, \lambda_s$ of the Weingarten map $S_{\eta}$ are called the \emph{principal curvatures} of $\Var X$ in the direction $\eta$. They measure how much~$\Var{X}$ curves at $X$ in the direction $\eta$. If $\lambda_i$ is a principal curvature of $\Var{X}$ at $X$ with associated unit-norm eigenvector $u_i$, then in the plane $P_i = \mathrm{span}(u_i, \eta)$ spanned by $u_i$ and $\eta$, the intersection of the manifold $\Var{X}$ with $P_i$ can be locally approximated to \textit{second} order at $X$ by a segment of an \textit{osculating circle} with center $X + \lambda_i^{-1} \eta$. This circular arc passes through $X$ with derivative $u_i \in \Tang{X}{\Var{X}}$; see \cref{fig_osculating}.

With this additional terminology, we obtain the following observation from \cref{eqn_Hess_distance}.

\begin{lemma}\label{norm_of_H}
Let $N:=A-X \in \Norm{X}{\Var{X}}$. Then,
\[
\Vert (H_{A,X})^{-1}\Vert_2 = \max_{1\leq i\leq s} (\vert 1 - \Vert N\Vert \lambda_i|)^{-1},
\]
where the norm is the spectral norm and the $\lambda_i$ are the principal curvatures.
\end{lemma}

\subsection{The second fundamental form as a tensor}\label{sec:SFF_3_tensor}
By multilinear algebra \cite{Greub1978}, we can represent the second fundamental form $\SFF_{X}$ from \cref{SFF} by a three-dimensional tensor in $(\Tang{X}{\Var{X}})^*\otimes (\Tang{X}{\Var{X}})^* \otimes \mathrm{N}_{X}\Var{X}$, where $(\,\cdot\,)^*$ denotes the dual. This tensor is symmetric in the first two factors; see \cite[equation (H)]{BV2020} or \cite{Petersen} for more details.

The dual space $(\Tang{X}{\Var{X}})^*$ is identified with $\Tang{X}{\Var{X}}$ via the standard Euclidean inner product $\langle\cdot,\cdot\rangle$ on $\R^n$ because $\mathcal X\subset \mathbb R^n$ is an embedded manifold inheriting the Riemannian structure from $\mathbb R^n$.
Let~$\vect{s}_i, \vect{t}_i \in \Tang{X}{\Var{X}}$ and $\vect{u}_i \in \Norm{X}{\Var{X}}$ be vectors so we can write
\begin{equation*}
\SFF_X = \sum_{i=1}^r \vect{s}_i \otimes \vect{t}_i \otimes \vect{u}_i
\end{equation*}
for some $r$; such an expression exists \cite{Greub1978}. The corresponding bilinear map is then
$
 \SFF_X( \vect{a}, \vect{b} ) = \SFF_X( \vect{b}, \vect{a} ) = \sum_{i=1}^r  \left( \langle \vect{s}_i, \vect{a} \rangle \cdot \langle \vect{t}_i, \vect{b} \rangle \right) \vect{u}_i.
$
Let $N \in \Norm{X}{\Var{X}}$ be a normal vector at $X$.
The contraction of $\SFF_X$ with $N$ along the third factor is defined by
\begin{equation*}
 N^T \cdot_3 \SFF_X
 = \sum_{i=1}^r \langle N, \vect{u}_i \rangle\, \vect{s}_i \otimes \vect{t}_i.
\end{equation*}
In the standard basis of $\R^n$, the tensor $\vect{s}_i \otimes \vect{t}_i$ would be represented by the rank-one matrix~$\vect{s}_i \vect{t}_i^T$, so that the foregoing equation can also be viewed naturally as a $\dim\Var{X} \times \dim\Var{X}$ matrix.
Comparing with the definition of the Weingarten map and \cref{eqn_Hess_distance}, we see the latter can also be expressed as
\begin{equation}\tag{H'}\label{H_new_form_general_setting}
H_{A,X} = \mathbf 1_{\Tang{X}{\Var X}} - N^T \cdot_3 \SFF_X.
\end{equation}

\section{Sensitivity of low-rank approximation}\label{sec_low_rank_approximation}
In low-rank matrix approximation the sensing operator $L$ in \cref{eqn_sensing_op} is given by $M$ being the identity on ${\Var{M}_r}$ and $b=0$.
Hence, $\Var{R}_r = \Var{S}_r$. It follows from \cref{def_CN} that the condition number of low-rank approximation is
\[
\kappa_\mathrm{approximation}(A,Y) = \|H_{A,Y}^{-1} \|_2,
\]
where $H_{A,Y}$ is the Riemannian Hessian of $d_A:\mathcal M_r \to \R,\, Y\mapsto \tfrac{1}{2}\Vert A-Y\Vert^2$ at $Y$.

Let $s:=\dim \Var M_r$ and $N=A-Y$. Since $Y$ minimizes the distance function $d_A$, we have $N\in\mathrm{N}_X \Var{R}_r$, i.e., $N$ is a normal vector of $\Var{R}_r$ at $Y$. By \cref{norm_of_H}, we have
\begin{equation}\label{cn_Mr}
\kappa_\mathrm{approximation}(A,Y) = \max_{1\leq i\leq s} (\vert 1 - \Vert N\Vert \lambda_i\vert)^{-1},
\end{equation}
where $\lambda_1,\ldots,\lambda_s$ are the principal curvatures of $\Var{R}_r$ at $Y$ and in direction $\eta:=\frac{N}{\Vert N\Vert}$.

The principal curvatures (of open submanifolds) of $\Var{M}_r \subset \R^{m \times n}$ can be derived from Amelunxen and B\"urgisser's Proposition 6.3 in \cite{AmBu2015} and were also stated by Feppon and Lermusiaux \cite[Theorem 24]{FL2018}. Let $A = \sum_{i=1}^{\min\{m,n\}} \sigma_i \vect{u}_i \vect{v}_i^T$ be an SVD of $A$ with the singular values $\sigma_1\geq \cdots\geq \sigma_{\min\{m,n\}}$ sorted decreasingly, and let the corresponding truncated rank-$r$ SVD be $Y =\sum_{i=1}^r\sigma_i \vect{u}_i \vect{v}_i^T \in \Var{R}_r$.
The principal curvatures at $Y$ in the normal direction $\eta = \frac{N}{\|N\|}$ are, on the one hand,
\begin{equation}\label{pc_of_Mr}
 c_{(b,i,j)} := \frac{(-1)^b}{\Vert N\Vert} \frac{\sigma_{{r+j}}}{\sigma_{i}}, \quad b=0,1;\; i = 1, \ldots, r; \text{ and } j=1,\ldots, \min\{m,n\}-r,
\end{equation}
and the other hand $r(m+n-r)-2r(\min\{m,n\}-r)$ principal curvatures are equal to zero \cite{AmBu2015,FL2018}. We can now prove \cref{main_thm1}.

\begin{proof}[Proof of \cref{main_thm1}]
We combine and \cref{cn_Mr} and \cref{pc_of_Mr} to get
\[
 \kappa_{\mathrm{approximation}}(A,Y)
 = \max_{\substack{1\leq i\leq r,\\ 1\leq j\leq \min\{m,n\}-r}} \Big(1 - \frac{\sigma_{{r+j}}}{\sigma_{i}}\Big)^{-1} =  \Big(1 - \frac{\sigma_{{r+1}}}{\sigma_{r}}\Big)^{-1}.
\]
This proves \cref{main_thm1}.
\end{proof}
\begin{remark}
In \cref{rmk1} we mentioned that the analysis of condition numbers also carries over the critical points. The critical points of the squared distance function $d_A$ for $A = \sum_{i=1}^{\min\{m,n\}} \sigma_i \vect{u}_i \vect{v}_i^T$ are all of the form $Y = \sum_{i\in I} \sigma_i \vect{u}_i \vect{v}_i^T$, where $I$ is a subset the indices with $|I|=r$. It can be shown that the condition number of low-rank approximation at such a critical point is $\max_{i\in I, j\not\in I} (1 - \tfrac{\sigma_{{j}}}{\sigma_{i}})^{-1}$.
\end{remark}

\Cref{main_thm1} shows that if there is a clear gap between the $r$th and $(r+1)$th singular value, then the best rank-$r$ approximation problem is well-conditioned. However, if $\sigma_r \approx \sigma_{r+1}$ then the problem is nearly ill-conditioned, by \cref{main_thm1}.

In the introduction we presented an informal example of an ill-conditioned low-rank approximation. For completeness, this example is formalized next.

\begin{example}\label{Example}
Let $0 < \epsilon \le 1$ be fixed and consider the matrix $A_0=\operatorname{diag}(1+\epsilon,1-\epsilon)$. The condition number at the unique best rank-$1$ approximation $Y_0=\operatorname{diag}(1+\epsilon,0)$ is
\(
\kappa_{\mathrm{approximation}}(A_0,Y_0) = \left( 1 - \frac{1-\epsilon}{1+\epsilon} \right)^{-1} = \frac{1}{2} (1 + \epsilon^{-1}).
\)
Let $\delta = \tau \epsilon$ with $0 \le \tau \le 1$. Consider the perturbed matrix $A_\delta = A_0 + \delta\left[\begin{smallmatrix}0 & 1 \\ 1 & 0 \end{smallmatrix}\right]$. Its eigendecomposition is
$$ A_\delta = V_\delta \Lambda_\delta V_\delta^{-1}$$
with
\begin{align*}
 \Lambda_\delta &= \diag(1 - \epsilon \sqrt{1 + \tau^2}), 1 + \epsilon\sqrt{1+\tau^2} )\quad \text{ and}\\
 V_\delta
 &=  \frac{1}{\tau} \begin{bmatrix}
1 - \sqrt{1 + \tau^2} & 1 + \sqrt{1+\tau^2} \\
 \tau & \tau
 \end{bmatrix}.
\end{align*}
Let $\zeta = \sqrt{1+\tau^2}$. The unique best rank-$1$ approximation of $A_\delta$ is then
 \begin{align*}
  Y_\delta
  &= \frac{1+\epsilon \zeta}{\tau^{2} + \left( 1 + \zeta \right)^2} \begin{bmatrix} 1+\zeta \\ \tau \end{bmatrix} \begin{bmatrix} 1+\zeta \\ \tau \end{bmatrix}^T
  = \frac{1+\epsilon\zeta}{\tau^{2} + \left( 1 + \zeta \right)^2}
  \begin{bmatrix}
    (1+\zeta)^2 & \tau(1+\zeta) \\
     \tau(1+\zeta) & \tau^2
     \end{bmatrix}.
 \end{align*}
Consequently,
\begin{align*}
 \|Y_\delta-Y_0\|_F^2
 &= \left(\frac{\tau(1+\epsilon\zeta)}{\tau^2 + (1+\zeta)^2}\right)^2 \left( 2 (1+\zeta)^2 + \tau^2 \right)+ \left( \frac{(1+\epsilon\zeta)(1+\zeta)^2}{\tau^2 + (1+\zeta)^2} - (1+\epsilon) \right)^2\\
 &= \frac{1}{2} (1+\epsilon)^2 \tau^2 + \frac{1}{8} (\epsilon^2 - 4 \epsilon -3)\tau^4 + \frac{1}{16} (-\epsilon^2 +6\epsilon + 5) \tau^6 + O(\tau^8),
\end{align*}
where the second equality is the Taylor series expansion around $\tau=0$.
Dividing by $\|A_\delta-A_0\|_F^2 = 2\delta^2 = 2 \tau^2 \epsilon^2$, we find
\[
 \frac{\|Y_\delta-Y_0\|_F^2 }{\|A_\delta-A_0\|_F^2}
= \frac{1}{4} (1+\epsilon^{-1})^2 + \frac{1}{16} (1 - 4\epsilon^{-1} -3\epsilon^{-2}) \tau^2 + O(\tau^{4}).
\]
The perturbation $\delta\left[\begin{smallmatrix}0&1\\1&0\end{smallmatrix}\right]$ of size $\sqrt{2\delta}$ that takes $A_0$ to $A_\delta$ moves the best rank-$1$ approximation from $Y_0$ to $Y_\delta$. The distance between these approximations relative to the distance between the matrices is approximately $\frac{1}{2} (1+\epsilon^{-1})$ plus higher-order terms in $\tau$. As $\epsilon$ was fixed and arbitrary, the low-rank approximation of $A_0$ can be made as ill-conditioned as wanted. Note that for $\epsilon=0$ the condition number tends to $\infty$, which in this case is caused by the occurrence of a positive-dimensional family of best rank-$1$ approximations of $\mathrm{diag}(1,1)$.

As a final observation, note that the limit for $\delta\to0$ of $\frac{\|Y_\delta - Y_0\|_F^2}{\|A_\delta - A_0\|_F^2}$ is precisely the square of the condition number $\kappa_\mathrm{approximation}(A_0,Y_0)$. That is, the perturbation $\delta\left[\begin{smallmatrix}0&1\\1&0\end{smallmatrix}\right]$ is exactly the worst direction of perturbation for $A_0$.
\end{example}

\section{Sensitivity of low-rank matrix recovery}\label{sec:recovery}
We continue with our discussion of low-rank recovery. Here, $L(Y) = M(Y) + b$ is a sufficiently general sensing operator for which we assume \cref{ass1} and \cref{prop_ass1} hold. The point $A\in \R^\ell$ is the input data to the recovery problem, $X=L(Y) \in\Var{S}_r$ is a sensed rank-$r$ matrix approximating $A$, and~$Y \in \Var{R}_r$ is the recoverable rank-$r$ matrix that projects to $X$.

Recall from \cref{def_CN} that the condition number of low-rank matrix recovery is
\[
\kappa_\mathrm{recovery}(A,Y) = \| (M|_{\mathrm T_Y\Var{R}_r})^{-1} H_{A,X}^{-1} \|_2,
\]
where $H_{A,X}$ is the Riemannian Hessian of the squared distance to the manifold $d_A : \Var{S}_r \to \R, X \mapsto \frac{1}{2} \| A - X \|^2$ at $X = L(Y)$.
This section derives a closed expression for $H_{A,X}$. Unfortunately, we are unable to derive a closed expression for $\kappa_\mathrm{recovery}(A,Y)$. Therefore, we will present a simple and efficient numerical linear algebra algorithm for evaluating $\kappa_\mathrm{recovery}(A,Y)$ in the next section.

Recall from \cref{H_new_form_general_setting} that the Riemannian Hessian can be expressed in terms of the second fundamental form.
Thus, our problem reduces to computing the latter. We can rely on the following lemma, which shows how curvature transforms under affine linear diffeomorphisms. While this is considered an elementary result in differential geometry, we could not locate a suitable reference, so a proof is included in the appendix for self-containedness.

\begin{lemma}
\label{lemma_L}
Consider Riemannian embedded submanifolds $\Var{U}\subset \R^N$ and $\Var{W} \subset \R^\ell$ both of dimension~$s$. Let $L : \R^N \to \R^\ell,\; Y\mapsto M(Y)+b$ be an affine linear map that restricts to a diffeomorphism from $\Var{U}$ to $\Var{W}$. For a fixed $Y\in \Var U$, let $(E_1,\ldots, E_s)$ be a basis of $\Tang{Y}{\Var{U}}$. For each $1\leq i\leq s$ let $F_i = M( E_i )$. Then, $(F_1,\ldots,F_s)$ is a basis of the tangent space $\Tang{X}{\Var{W}}$ at $X=L(Y)$, and we have
 \[
\SFF_{X}(F_i, F_j) = \mathrm{P}_{\mathrm{N}_{X} \Var{W}}\left(  M( \SFF_Y(E_i, E_j) ) \right).
 \]
\end{lemma}

This lemma shifts our problem to computing the second fundamental form of (recoverable) rank-$r$ matrices $\Var{R}_r \subset \Var{M}_r$. The latter was computed in \cite[Section 4.5]{AMT2013} and \cite[Proposition 22]{FL2018}. In the next subsection we will evaluate the latter at an orthonormal basis, so a succinct matrix representation is obtained.

\subsection{Second fundamental form of rank-$r$ matrices}\label{SFF_rank_r}
Let $Y = U \Sigma V^T$ be a compact SVD of $Y \in \Var{R}_r$, such that $\Sigma$ is the diagonal matrix with entries the singular values $\sigma_1\geq \cdots\geq \sigma_r$.
Then, by the fact that $\Var{R}_r$ is an open submanifold of $\Var{M}_r$ and \cite{HM1994}, the tangent and normal spaces to $\Var{R}_r$ at $Y$ are
\begin{equation}\label{tangent_space_M_r}
\Tang{Y}{\Var{R}_r} = (U \otimes V) \oplus (U^\perp \otimes V) \oplus (U \otimes V^\perp)
\quad\text{ and }\quad
\Norm{Y}{\Var{R}_r} = U^\perp \otimes V^\perp,
\end{equation}
where $U$ and $V$ are conveniently identified with their column spans, $\oplus$ denotes the direct sum of (orthogonal) linear subspaces, and $(\cdot)^\perp$ denotes the orthogonal complement of a subspace. Let the columns of $U$ be $u_1, \ldots, u_r$, and let $v_1, \ldots, v_r$ be the columns of~$V$. Let $u_{r+1}, \ldots, u_m$ be an orthonormal basis of $U^\perp$, and $v_{r+1},\ldots,v_{n}$ one for $V^\perp$. Then,
\begin{alignat*}{2}
 U &\otimes V &&= \operatorname{span}( u_i v_j^T \mid 1 \le i, j \le r ),\\
 U^\perp &\otimes V &&= \operatorname{span}( u_i v_j^T \mid 1 \le j \le r < i \le m ), \\
 U &\otimes V^\perp &&= \operatorname{span}( u_i v_j^T \mid 1 \le i \le r < j \le n ),\\
 U^\perp &\otimes V^\perp &&= \operatorname{span}( u_i v_j^T \mid r < i \le m, r < j\le n ).
\end{alignat*}
For brevity we define $E_{ij} = u_i v_j^T$ for all $i$ and $j$. We also define the rank-$1$ matrices
\[
\phi_{ij} = \begin{cases}
0 & \text{if } E_{ij} \in U\otimes V, \\
0 & \text{if } E_{ij} \in U\otimes V^\perp,\\
\sigma_j^{-1} u_i e_j^T & \text{if } E_{ij} \in U^\perp \otimes V,
\end{cases}\quad\text{and}\quad
\psi_{ij} =
\begin{cases}
v_j e_i^T & \text{if } E_{ij} \in U\otimes V, \\
v_j e_i^T & \text{if } E_{ij} \in U\otimes V^\perp,\\
0 & \text{if } E_{ij} \in U^\perp \otimes V.
\end{cases}
\]
Then, we have the unique decomposition $E_{ij} = U \psi_{ij}^T + \phi_{ij} (V \Sigma)^T$.
We obtain the following formula from \cite[Proposition 22]{FL2018}:
\[
\SFF_Y( E_{i,j}, E_{k,l} )
= \mathrm{P}_{\Norm{Y}{\Var{M}_r}}( \phi_{ij}\, \psi_{kl}^T +  \phi_{kl} \psi_{ij}^T).
\]
It can be verified by direct computation that the expression simplifies to
\begin{align} \label{eqn_sff_compact}
 \SFF_Y( E_{ij}, E_{kl} ) =
 \begin{cases}
\sigma_j^{-1} \delta_{kj} E_{il}  & \text{if } E_{ij} \in U^\perp \otimes V \text{ and } E_{kl} \in U\otimes V^\perp, \\
\sigma_{l}^{-1} \delta_{il} E_{kj}  & \text{if } E_{ij} \in U \otimes V^\perp \text{ and } E_{kl} \in U^\perp\otimes V, \\
0 & \text{otherwise}.
 \end{cases}
\end{align}
Herein $\delta_{ab}$ is the Kronecker delta.
The fact that $\SFF_Y$ restricted to $U\otimes V$ is zero is actually a priori clear from geometric considerations: the second fundamental form measures the curvature of $\Var{R}_r$ inside of $\R^{m\times n}$ and if $Y = U \Sigma V^T$, there is a whole linear space contained in $\Var{R}_r$ passing through $Y$, namely $U \otimes V$. The part of the second fundamental form arising from contravariant differentiation of the basis vectors of $U \otimes V$ thus vanishes completely.

Since the $E_{ij}$ form an orthonormal basis, it can be deduced from \cref{eqn_sff_compact} that the following is the second fundamental form of $\Var{M}_r$ at $Y$ viewed as element of the tensor space $(\Tang{Y}{\Var{R}_r})^* \otimes (\Tang{Y}{\Var{R}_r})^* \otimes \Norm{Y}{\Var{R}_r}$:
\begin{align}
 \nonumber
 \SFF_Y &= \sum_{i=r+1}^m \sum_{j=1}^r \sum_{k=1}^r \sum_{l=r+1}^n E_{ij}^T \otimes E_{kl}^T \otimes \Big( \frac{1}{\sigma_j} \delta_{kj} E_{il} \Big)  \\ \nonumber
 &\hspace{4.25cm} + \sum_{i=1}^r \sum_{j=r+1}^n \sum_{k=r+1}^m \sum_{l=1}^r E_{ij}^T \otimes E_{kl}^T \otimes \Big( \frac{1}{\sigma_l} \delta_{il} E_{kj} \Big),\\
 \nonumber
 &= \sum_{i=r+1}^m \sum_{l=r+1}^n \sum_{k=1}^r \frac{1}{\sigma_k} E_{ik}^T \otimes E_{kl}^T \otimes E_{il} + \sum_{k=r+1}^m  \sum_{j=r+1}^n \sum_{l=1}^r \frac{1}{\sigma_l} E_{lj}^T \otimes E_{kl}^T \otimes E_{kj}, \\
\label{eqn_stuff_complete} &= \sum_{i=r+1}^m \sum_{j = r+1}^n \sum_{k=1}^r \frac{1}{\sigma_k} (E_{ik}^T \otimes E_{kj}^T + E_{kj}^T \otimes E_{ik}^T) \otimes E_{ij};
\end{align}
see also the discussion in \cref{sec:SFF_3_tensor}.
Note that since our Riemannian metric is the standard Euclidean inner product $(A,B)\mapsto \mathrm{Trace}(A^TB)$ on $\R^{m \times n}$, dualization consists of transposition.

\subsection{Second fundamental form of sensed rank-$r$ matrices}\label{SFF_incomplete_rank_r}
We can now compute the second fundamental form of the sensed manifold $\Var{S}$.

Let us denote $F_{ij} = M(E_{ij})$. These are the images of the basis vectors of $\Tang{Y}{\Var{R}_r}$ under the derivative $\deriv{L}{Y}( \dot{Y} ) = M(\dot{Y})$. As before, let $X=L(Y) \in\Var{S}_r$.
We conclude from \cref{lemma_L,eqn_stuff_complete} that the second fundamental for $\SFF_X$, viewed as an element of $(\Tang{X}{\Var{S}_r})^* \otimes (\Tang{X}{\Var{S}_r})^* \otimes \Norm{X}{\Var{S}_r}$, is
\begin{equation}\label{SFF_X}
 \SFF_{X} = \sum_{i=r+1}^m \sum_{j = r+1}^n \sum_{k=1}^r \frac{1}{\sigma_k} \left( F_{ik}^\dagger \otimes F_{kj}^\dagger + F_{kj}^\dagger \otimes F_{ik}^\dagger \right) \otimes \mathrm{P}_{\Norm{X}{\Var{W}_r}}(M(E_{ij})),
\end{equation}
where $F_{ij}^\dagger$ is the dual basis vector of $F_{ij}$; that is, $\langle F_{ij}^\dagger, F_{kl} \rangle = \delta_{ik} \delta_{jl}$.\footnote{It is customary to denote the dual basis by $F_{ij}^*$. This dual basis is often defined with respect to an orthonormal basis in $\mathbb R^\ell$. We chose $\dagger$ to emphasize that taking the dual of $F_{ij}$ in this way does not result in the dual basis vector $F_{ij}^\dagger$.}
Recall from~\cref{H_new_form_general_setting} that the formula for the Riemannian Hessian is $H_{A,X} = \mathbf 1_{\Tang{X}{\Var X}} - N^T \cdot_3 \SFF_X$, where $N=A - X$. The second term in this formula is thus
\begin{equation}\label{SFF_X_2}
N^T \cdot_3 \SFF_X = \sum_{i=r+1}^m \sum_{j = r+1}^n \sum_{k=1}^r \frac{1}{\sigma_k} \langle N, M(E_{ij}) \rangle \, \left( F_{ik}^\dagger \otimes (F_{kj}^\dagger)^T + F_{kj}^\dagger \otimes (F_{ik}^\dagger)^T \right).
\end{equation}
As discussed in \cref{sec:SFF_3_tensor}, this expression can be represented naturally by a matrix in $(\Tang{X}{\Var{S}_r})^* \otimes \Tang{X}{\Var{S}_r}$. The transposition on the right-hand side originated from taking duals as $N^T \cdot_3 \SFF_X$ can be seen as a bilinear map $\Tang{X}{\Var{S}_r} \times \Tang{X}{\Var{S}_r} \to \R$.

Equation \cref{SFF_X_2} specifies in abstract terms the contraction of the second fundental form by $N$. We can use this to compute the condition number $\kappa_\mathrm{recovery}(A,Y)$. Indeed, $\kappa_\mathrm{recovery}(A,Y)$ is the spectral norm of the inverse of~$H_{A,X} \circ M$. Unfortunately, we were not able to determine a handy expression for the inverse of $H_{A,X}$. For this reason, we explain how the condition number can be computed using standard linear algebra software in the next subsection.

\section{An algorithm for computing the condition number}\label{alg_matrix_recovery}

The spectral norm in the definition \cref{def_CN} of $\kappa_\mathrm{recovery}(A,Y)$ can be computed efficiently in coordinates if we choose orthonormal bases for respectively the codomain and domain of the operator $H_{A,X}M$. In such bases, the spectral norm coincides with the $2$-norm of the coordinate matrix by classic linear algebra. The inverse of the smallest singular value of this matrix representation of $H_{A,X}M$ is then the condition number $\kappa_\mathrm{recovery}(A,Y)$.

\subsection{Determining the dual basis} First, we express the dual basis $F_{ij}^\dagger \in (\Tang{X}{\Var{S}_r})^*$ in the standard basis of $(\R^\ell)^*$.
We assume that $M$ computes the coordinates in the standard basis $(e_1, \ldots, e_\ell)$ of $\R^\ell$. Consequently, the basis vectors $F_{ij}$ are given in these coordinates by
\(
 F_{ij} = M( E_{ij} ).
\)
Let $F = [F_{ij}]$ be the $\ell \times s$ matrix formed by placing the $F_{ij}$'s as column vectors, where $s = \dim\Var{R}_r = (m+n-r)r$.

The dual basis of $F$, expressed in coordinates with respect to the standard basis $(e_1^T, \ldots, e_\ell^T)$ of $(\R^\ell)^*$, is then given by the rows of the Moore--Penrose pseudoinverse of $F$; indeed, $F^\dagger F = I_{s}$ so the rows $F_{ij}^\dagger$ are the dual basis vectors.

\subsection{Matrix representation of the Weingarten map}

From \cref{SFF_X_2}, we can now conclude that the matrix of the Weingarten map $S_N=N^T \cdot_3 \SFF_X$ relative to the standard basis $(e_1,\ldots,e_\ell)$ of $\R^\ell$ and $(e_1^T,\ldots,e_\ell^T)$ of $(\R^\ell)^*$ is
\[
S_N =
(F^\dagger)^T \begin{bmatrix}
 0_{r^2} &  &  \\[0.2em]
 & 0_{(m-r)r} & V \\[0.3em]
  & V^T & 0_{(n-r)r}
 \end{bmatrix}
 F^\dagger,
\]
where $0_{a}$ denotes an $a \times a$ matrix of zeros, $V$ is defined in the next paragraph, and all non-displayed entries are zero. Consequently, $S_N$ is a square matrix of size $r^2 + (m-r)r + (n-r)r = (m+n-r)r = s= \dim \mathcal M_r$.

We see from \cref{SFF_X_2} that the foregoing matrix $V \in \R^{(m-r)r \times (n-r)r}$ is indexed by a multi-index $(I,J) = ((ij),(kl))$ with $1 \le j \le r < i \le m$ and $1 \le n \le r < l \le n$. Its entries are:
\begin{align} \label{eqn_compute_v}
 V = \begin{bmatrix} \frac{\delta_{jk}}{\sigma_j} \langle N, M(E_{il}) \rangle \end{bmatrix}_{(ij),(kl)};
\end{align}
compare this with the start of the equations that led to \cref{eqn_stuff_complete}.

\begin{remark}
  Note that for fixed $(i, l)$, the submatrix of $V$ formed by $1 \le j, k \le r$ is a multiple of the identity. After a suitable symmetric permutation of rows and columns, we thus can write $V  = \Sigma^{-1} \otimes Z$, where $\Sigma = \operatorname{diag}(\sigma_1, \ldots, \sigma_r)$ and $Z = [\langle N, M(u_{i} v_l^T)\rangle ]_{\substack{r+1 \le i \le m,\\ r+1 \le l \le n}}$.
  This observation can be exploited to further simplify computations with $V$. Doing this implies a particular ordering of the basis vectors in~$(F^\dagger)^T$, which needs to be respected when computing $R$ below. As the computational gains associated with this observation do not lead to an improvement of the asymptotic running time, we decided not to exploit it in the discussion below.
\end{remark}

Consider the $QR$ factorization $F = QR$. As the columns $F \in \R^{\ell \times s}$ form a basis (recall that $\ell \ge s$), $R$ is an $s \times s$ invertible matrix. It follows that $F^\dagger = R^{-1} Q^T$. We can thus write
\[
 S_N = Q R^{-T} \begin{bmatrix}
 0 & 0 & 0 \\
 0 & 0 & V \\
 0 & V^T & 0
 \end{bmatrix}
 R^{-1} Q^T.
\]
Partitioning $R$ conformally with the block structure of the middle matrix, we have
\[
 R =
 \begin{bmatrix}
 R_{11} & R_{12} & R_{13} \\
 0 & R_{22} & R_{23} \\
 0 & 0 & R_{33}
 \end{bmatrix}, \quad\text{and}\quad
R^{-1} =
\begin{bmatrix}
 R_{11}^{-1} & T & T' \\
 0 & R_{22}^{-1} & -R_{22}^{-1} R_{23} R_{33}^{-1} \\
 0 & 0 & R_{33}^{-1}
 \end{bmatrix},
\]
where $T, T'$ are unspecified matrices. Consequently, the foregoing expression of $S_N$ can be simplified to
\[
S_N = Q \begin{bmatrix}
    0 & 0 & 0 \\
    0 & 0 & R_{22}^{-T} V R_{33}^{-1} \\
    0 & R_{33}^{-T} V^T R_{22}^{-1} & - \mathrm{Sym}(R_{33}^{-T} V^T R_{22}^{-1} R_{23} R_{33}^{-1} )
   \end{bmatrix}
 Q^T,
\]
where $\mathrm{Sym}(Z) = Z + Z^T$ symmetrizes its input.

\subsection{Computing the condition number}
Recall that the columns of $Q$ form an orthonormal basis of $\Tang{X}{\Var{S}_r}$. Therefore, the identity $\mathbf{1}_{\Tang{X}{\mathcal S_r}}$ is represented with respect to the standard basis on $\R^\ell$ (and its dual) as $Q Q^T$. By definition, $F = QR$ is the change of basis matrix $M : \Tang{Y}{\Var{R}_r} \to \Tang{X}{\Var{S}_r}$ from $E_{ij}$ to $F_{ij}$ represented with respect to the orthonormal basis $E_{ij}^T \in (\R^{m\times n})^*$ and the standard basis on $\R^\ell$. Putting all of this together, and using
\cref{eqn_Hess_distance}, we find that
\begin{align} \label{eqn_compute_tn}
 H_{A,X} M =
 Q
\begin{bmatrix}
  \mathbf 1_{r^2} & \\
  & \mathbf 1_{(m-r)r} & - R_{22}^{-T} V R_{33}^{-1} \\
  & -R_{33}^{-T} V^T R_{22}^{-1} & \mathbf{1}_{(n-r)r} + \mathrm{Sym}(R_{33}^{-T} V^T R_{22}^{-1} R_{23} R_{33}^{-1} )
 \end{bmatrix}
  R,
\end{align}
where $\mathbf 1_a$ is the $a\times a$ identity matrix. Let us write $T_N$ for the matrix in the middle in \cref{eqn_compute_tn}, so that $H_{A,X}M = QT_NR$.
This is a matrix representation relative to the standard orthogonal basis of $\R^\ell$ and the orthonormal basis $E_{ij}^T$ on~$(\R^{m\times n})^*.$ 
It follows that
\[
 \kappa_{\mathrm{recovery}}(A,Y)
 = \frac{1}{\sigma_{s}( T_N R )},
\]
where $\sigma_s(T_N R)$ is the smallest singular value of the $s \times s$ matrix $T_N R$ and where, as before, $s = \dim \Var{M}_r = (m+n-r)r$. We can ignore $Q$ because it has orthonormal columns.

\subsection{The algorithm}
We can now put all components together. We assume that a rank-$r$ matrix $Y = U \Sigma V^T \in \R^{m \times n}$ is given (factored or not). We assume without loss of generality that $m \ge n$. We are also given $X = L(Y)$ and $N = A - X$ lies (approximately) in the normal space $\Norm{X}{\Var{S}_r}$. As before, $s = \dim \Var{R}_r$.

The numerical algorithm we propose for the condition number $\kappa_\mathrm{recovery}(A,Y)$ proceeds as follows:
\begin{enumerate}
 \item[S1.] Compute orthonormal bases $U^\perp \in \R^{m \times m-r}$ and $V^\perp \in \R^{n \times n-r}$ for the orthogonal complements of $U$ and $V$ respectively via a full SVD of $Y$.
 \item[S2.] Construct the $\ell \times s$ change of basis matrix $F = [M(E_{ij})]_{i \le r \text{ or } j \le r}$ as well as the $\ell \times mn - s$ matrix $G = [M(E_{ij})]_{i, j > r}$.
 \item[S3.] Compute the QR decomposition $F = Q R$.
 \item[S4.] Ensure that $N = A - X$ is numerically orthogonal to the tangent space $\Tang{X}{\Var{S}_r}$ by computing $N \leftarrow N - Q (Q^T N)$ twice.
 \item[S5.] Construct the matrix $V$ by the formula \cref{eqn_compute_v} and the precomputed $G$.
 \item[S6.] Compute the matrix $T_N$ following \cref{eqn_compute_tn}, and then compute  $Z = T_N R$.
 \item[S7.] Compute the smallest singular value $\sigma_s(Z)$ of $Z$.
   \item[S7.] Output $\kappa_\mathrm{recovery}(A,Y) = \sigma_s(Z)^{-1}$.
\end{enumerate}

The cost of computing the condition number with the foregoing algorithm depends on the cost of applying the linear part $M$ of the sensing operator to the basis vectors $E_{ij} = u_i v_j^T$. Let us denote the maximal cost by $C_M$. The cost is
\begin{align*}
 &\underbrace{ m^3}_{S1.} + \underbrace{mn C_M}_{S2.} + \underbrace{\ell s^2}_{S3.} +  \underbrace{\ell s}_{S4.} + \underbrace{\ell (mn - s) r}_{S5.} + \underbrace{s^3}_{S6.} + \underbrace{s^3}_{S7.} \\[0.3em]
 = \; &\mathcal{O}( mn C_M + \ell s^2 ),
\end{align*}
where in the last step we used $\ell > s = (m+n-r)r$ and $(mn-s)r < s^2$.
For practical sensing operators with $\ell = \phi s$ and $\phi > 1$ a small constant, this usually means the cost is dominated by the cost for computing the QR-factorization of the change-of-basis matrix $F$.

A general sensing operator $L : \R^{m\times n} \to \R^\ell$ has cost $C_M = mn\ell$. As we have $\ell > (m+n-r)r$, this implies the overall cost for computing the condition number would be a rather impressive $m^3 n^2 r$. Fortunately, many sensing operator are structured. Consider, for example, the structured sensing operator
 \begin{align} \label{eqn_structured_sensing}
 L(X) = \operatorname{diag}(B^T X C) + \vect{b} = (B \odot C)^T \vecc{X} + \vect{b},
 \end{align}
which is defined by the $m \times \ell$ matrix $B$, the $n \times \ell$ matrix $C$ and a vector $\vect{b} \in \R^\ell$. In the foregoing, $\odot$ is the (columnwise) Khatri--Rao product of its arguments.
The derivative of $L$ is the map $\dot{X} \mapsto \diag(B^T \dot{X} C)$. Hence, if $\dot{X} = u_i v_j^T$, we see that the derivative can be applied effectively by $(B^T u_i) \circledast (v_j^T C)$, where $\circledast$ is the Hadamard or elementwise product. The computational complexity is only $\ell(m+n+1)$ in this case. With such a sensing operator, the condition number can be computed in~$\mathcal{O}(s^3)$ operations. In conclusion, we proved the next result.

\begin{proposition}\label{prop_complexity}
 Let the sampling operator be as in \cref{eqn_structured_sensing} and $\ell = \phi s$. Then, the condition number $\kappa_\mathrm{recovery}(A,Y)$ where $A \in \R^\ell$ and $Y \in \Var{R}_r \subset \R^{m\times n}$ can be computed in $\mathcal{O}(\phi s^3)$ operations, where $s = \dim \Var{M}_r = (m+n-r)r$.
\end{proposition}

This complexity is cubic in the problem size $s = \dim\Var{M}_r$. This means that computing the solution's condition number is as expensive as one step of a Riemannian Newton method for solving the recovery problem.

An example of such a structured sensing operator appears in the Netflix problem from the introduction. Herein, the sensing operator $L$ selects $\ell$ coordinates $(i_k,j_k)$ of $\R^{m\times n}$ and the other elements are unknown. This can be expressed as in \cref{eqn_structured_sensing} by taking $\vect{b}=0$, $B = [e_{i_k}]_k$ and $C = [e_{j_k}]_k$.

\section{Numerical experiment}\label{sec:experiments}

We present an experiment to study the condition number of low-rank matrix recovery. It was performed on a computer running Ubuntu 18.04.5 LTS, comprising a quad-core Intel Core i7-4770K CPU (3.5GHz clockspeed) and 32GB main memory. Our Julia implementation including experiments is available from the repository \url{https://gitlab.kuleuven.be/u0072863/MatrixRecoverySensitivity}.

We investigate the sensitivity of low-rank matrix recovery where the sensing operator is a random low-rank sensing operator. The $i$th measurement of the sensing operator $L : \R^{m \times n} \to \R^\ell$ performs $L_i(Y) = \mathbf{v}_i^T Y \mathbf{w}_i$ where $\mathbf{v}_i \in \R^m$ and $\vect{w}_i \in \R^n$ are vectors whose elements are drawn i.i.d.~from a standard Gauss distribution. We investigate the influence of the number of measurements
$$\ell = \varphi \dim \Var{M}_r = \varphi (m+n-r)r.$$
Here, $\varphi$ is the \textit{oversampling rate}: $\varphi=1$ is expected to suffice for finite recoverability; see~\cite[Section 3]{BGMV2021}.

We also investigate the influence the relative distance $t$ of the input matrix
\[
A_t = X +   t\frac{\|X\|}{\Vert N\Vert}\cdot N
\]
from the sensed input manifold $\Var{S}_r$. Herein, $X=L(Y) \in \R^\ell$ is the image under~$L$ of a randomly chosen rank-$r$ matrix $Y=AB^T$ with $A$ and $B$ random Gaussian matrices, and $N$ is a random unit-norm normal vector at $Y$. The normal vector $N$ is chosen as follows: we sample $\eta$ as a random Gaussian vector in $\mathbb R^\ell$ and orthogonally project onto the normal space so that $N=\mathrm{P}_{\mathrm{N}_{X}\Var S_r}(\eta)$.

\begin{figure}[t]
\begin{center}
\includegraphics[height=6.625cm]{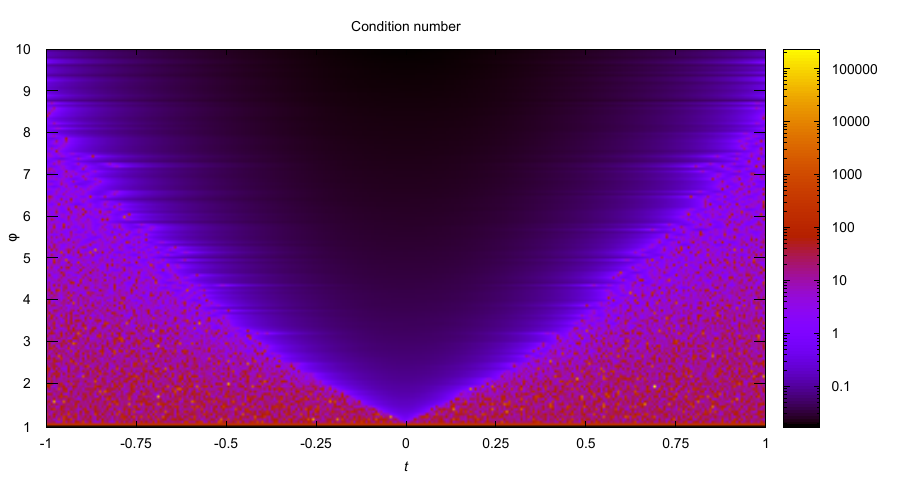}
\end{center}
\caption{The base-$10$ logarithm of the condition number $\kappa_\textrm{recovery}$ for various combinations of the oversampling factor $\phi$ and the distance signed $t$.}
\label{exp_LR3}
\end{figure}

In our experiment, we took $(m,n,r)=(50,40,10)$. For $t$ we took $299$ linearly spaced samples between $-1$ and $1$, and $\varphi$ $150$ linearly spread samples between $1$ and $10$ were chosen. The corresponding number of measurements $\ell$ was the integer part of $\varphi \dim \Var{M}_r$. Note that all $\ell_{\max} = 10\dim\Var{M}_r$ vectors $\vect{v}_i$ and $\vect{w}_i$ are generated beforehand and we always use the first $\ell$ measurements for a particular $\varphi$. We are thus only adding measurements as $\varphi$ is increased.

The base-$10$ logarithm of the condition number of low-rank recovery at $(A_t,Y)$ is visualized in \cref{exp_LR3}. By considering a vertical column in the figure, we can see the effect of adding additional measurements on the condition number. A phase transition can be made out in the figure. The very dark area in the figure corresponds to the cases where $Y$ is a local minimizer of the distance to $A_t$. In the purple--red area, on the other hand, $Y$ is no longer a local minimizer and the condition numbers can be significantly higher; anywhere from around $10$ to $\infty$.

The key feature \cref{exp_LR3} demonstrates is that some amount of oversampling $\varphi>1$ in the measurements is necessary for very well-conditioned local minimizers, i.e., $\kappa_\text{recovery}(A_t,Y) \le 1$, especially when the input matrix~$A_t$ does not lie on the manifold of sensed matrices $\Var{S}_r$, i.e., $|t|>0$. In many practical applications this is true, such as in the Netflix problem, because it is only \textit{assumed} that the output $Y$ can be well-approximated by a low-rank matrix on $\Var{R}_r$. Consequently, the sensed matrix is not expected to lie on the sensed manifold $\Var{S}_r = L(\Var{R}_r)$ either.

\appendix

\section{Proof of \cref{lemma_L}.}

We restate \cref{lemma_L} in terms of vector fields, as required by the proof, and then we prove it.

\smallskip

\noindent{\bf \cref{lemma_L}} \; (The second fundamental form under affine linear diffeomorphisms)\\
\noindent{\it
  Consider Riemannian embedded submanifolds $\Var{U}\subset \R^N$ and $\Var{W} \subset \R^\ell$ both of dimension~$s$. Let $L : \R^N \to \R^\ell, Y\mapsto M(Y)+b$ be an affine linear map that restricts to a diffeomorphism from $\Var{U}$ to $\Var{W}$. For a fixed $Y\in \Var U$ let $\Var E =(E_1,\ldots, E_s)$ be a local smooth frame of $\Var U$ in the neighborhood of $Y$. For each $1\leq i\leq s$ let $F_i$ be the vector field on $\Var W$ that is $L|_{\Var U}$-related to $E_i$. Then $\Var F = (F_1,\ldots,F_s)$ is a smooth frame on~$\Var W$ in the neighborhood of $X=L(Y)$, and we have
   \[
  \SFF_{X}(F_i, F_j) = \mathrm{P}_{\mathrm{N}_{X} \Var{W}}\left(  M( \SFF_Y(E_i, E_j) ) \right).
   \]
   Herein, $\SFF_Y \in \Tang{Y}{\Var{U}} \otimes \Tang{Y}{\Var{U}} \otimes \Norm{Y}{\Var{U}}$ should be viewed as an element of the tensor space $(\Tang{Y}{\Var{U}})^* \otimes (\Tang{Y}{\Var{U}})^* \otimes \Norm{Y}{\Var{U}}$
   by dualization relative to the Riemannian metric, and likewise for $\SFF_{X}$.}
\smallskip
\begin{proof}
Our assumption of $L|_{\Var U}:\Var U \to \Var W$ being a diffeomorphism implies that $\Var F$ is a smooth frame. Let $1\leq i,j\leq s$, and $\varepsilon_i(t) \subset \Var{U}$ be the integral curve of $E_i$ starting at $Y$ (see \cite[Chapter 9]{Lee2013}). Let $\phi_i(t) = L(\varepsilon_i(t))$ be the corresponding integral curve on~$\Var{W}$ at $X=L(Y)$. Since $F_j$ is~$L|_{\Var U}$-related to $E_j$, we have
\[
 F_j|_{\phi_i(t)}
 = (\deriv{L|_{\Var U}}{\varepsilon_i(t)}) (E_j|_{\varepsilon_i(t)}).
\]
where $\deriv{L|_{\Var U}}{\varepsilon_i(t)}$ is the derivative $\deriv{L}{\varepsilon_i(t)}: \R^N\to\R^\ell$ of $L$ at $\varepsilon_i(t)$ restricted to the tangent space $\mathrm{T}_Y\Var U\subset \R^N$.
On the other hand, $\deriv{L}{\varepsilon_i(t)}\ = M $, so that
\begin{equation}\label{F_related}
 F_j|_{\phi_i(t)}
 = M \, E_j|_{\varepsilon_i(t)}.
\end{equation}
The fact that the derivative of $L$ is constant is the key part in the proof.
Interpreting~$F_{j}|_{\phi_i(t)}$ as a smooth curve in $\Tang{\phi_i(t)}{\R^\ell} \simeq \R^\ell$ and $E_j|_{\varepsilon_i(t)}$ as a smooth curve in~$\Tang{\varepsilon_i(t)}{\R^N} \simeq \R^N$, we can take the usual derivatives at $t=0$ on both sides of \cref{F_related}:
\[
 \frac{\mathrm{d}}{\mathrm{d} t} F_{j}|_{\phi_i(t)}
 = M\, \frac{\mathrm{d}}{\mathrm{d} t} E_j|_{\varepsilon_i(t)}.
\]
Recall that $Y=\varepsilon_i(0)$.
We can decompose the right hand side into tangent and normal part at $Y\in\mathcal U$, so that
$$
\frac{\mathrm{d}}{\mathrm{d} t} F_{j}|_{\phi_i(t)}= M\, \left( \mathrm{P}_{\Tang{Y}{\Var{U}}} \frac{\mathrm{d}}{\mathrm{d} t} E_j|_{\varepsilon_i(t)} \right) \oplus M\, \left( \mathrm{P}_{\mathrm{N}_{x}\Var{M}} \frac{\mathrm{d}}{\mathrm{d} t} E_j|_{\varepsilon_i(t)} \right).
$$
Observe that that $M(\Tang{Y}{\Var{U}}) = \Tang{X}{\Var{W}}$, where $X=L(Y)$. Projecting both sides to the normal space of $\Var{W}$ at $X$ yields
\[
\mathrm{P}_{\mathrm{N}_{X} \Var{W}} \left(  \frac{\mathrm{d}}{\mathrm{d} t} F_{j}|_{\phi_i(t)}
 \right)
 = \mathrm{P}_{\mathrm{N}_{X} \Var{W}} \left( M\, \mathrm{P}_{\mathrm{N}_{Y}\Var{U}} \left( \frac{\mathrm{d}}{\mathrm{d} t} E_j|_{\varepsilon_i(t)} \right) \right).
\]
The claim follows by applying the Gauss formula for curves \cite[Lemma 8.5]{Lee1997} on both sides.
\end{proof}

\bibliographystyle{amsplain}
\bibliography{BV6}
\end{document}